\documentclass[journal]{ieeeconf}
\usepackage{ifpdf}
\usepackage{cite}
   \usepackage[pdftex]{graphicx}
   \usepackage{comment}
\usepackage{amsmath}
\usepackage{amssymb}
\usepackage[algoruled]{algorithm2e}
\usepackage{array}
  \usepackage[caption=false,font=footnotesize]{subfig}
\usepackage{url}
\usepackage{multicol}
\usepackage{booktabs}

\usepackage{color}
\usepackage{cite}

\begin{document}
\title{Koopman Operator Theory for Nonlinear Dynamic Modeling using Dynamic Mode Decomposition}

% author names and affiliations (conference format)

\author{Gregory~Snyder
	and~Zhuoyuan~Song% <-this % stops a space        
\thanks{*This work was supported by the U.S.~National Science Foundation under awards CISE/IIS-2024928 and OIA-2032522.}
\thanks{G.~Snyder and Z.~Song are with the Department of Mechanical Engineering, University of Hawai`i at M\={a}noa, Honolulu,
HI, 96822 USA. E-mails: \{snyderg, zsong\}@hawaii.edu.}% <-this % stops a space
}

\maketitle


\begin{abstract}

    % With advancements in computing power and the rising popularity of using data driven methods for system identification, the Koopman operator and techniques to approximate it have 
	The Koopman operator is a linear operator that describes the evolution of scalar observables (i.e., measurement functions of the states) in an infinite-dimensional Hilbert space.
	This operator theoretic point of view lifts the dynamics of a finite-dimensional nonlinear system to an infinite-dimensional function space where the evolution of the original system becomes linear.
	In this paper, we provide a brief summary of the Koopman operator theorem for nonlinear dynamics modeling and focus on analyzing several data-driven implementations using dynamical mode decomposition (DMD) for autonomous and controlled canonical problems. 
	We apply the extended dynamic mode decomposition (EDMD) to identify the leading Koopman eigenfunctions and approximate a finite-dimensional representation of the discovered linear dynamics. 
	This allows us to apply linear control approaches towards nonlinear systems without  linearization approximations around fixed points. 
	We can then examine the fidelity of using a linear controller based on a Koopman operator approximated system on under-actuated systems with basic maneuvers.
	We demonstrate the effectiveness of this theory through numerical simulation on two classic dynamical systems are used to show DMD methods of evaluating and approximating the Koopman operator and its effectiveness at linearizing these systems.
\end{abstract}

\iffalse
\fi

% \fi
\IEEEpeerreviewmaketitle

%\section{Introduction}

\section{Introduction}
%\IEEEPARstart{L}{ocalization} 

In this work we explore the application of the Koopman operator to find finite-dimensional linear representation of nonlinear dynamical systems. The Koopman operator is an infinite-dimensional operator that evolves the functions that evolve the state of a dynamical system linearly. In this work, the Koopman operator is computed using data-driven methods such as the dynamic mode decomposition (DMD) algorithm. DMD uses linear measurements of the states of a dynamical system to find the dominant modes of the underlying system from data to reconstruct a linear representation of the evolution of the system. 

The theory behind the Koopman operator was first postulated as an infinite-dimensional linear operator for uncontrolled systems in the seminal work by Koopman and Neumann\cite{Koopman1932}, which recently has sparked an reemergence of interest in this method of modeling systems in the past two decades. 
More recently, methods of finding the Koopman operator have been looked into by using discrete data-driven methods for fluid dynamics~\cite{rowley2009} through dimensional reduction algorithms developed by Schmid et al.~\cite{schmid_2010}.

Brunton et al.~explored the relationship between the Koopman operator and explored multiple observable functions to form a Koopman subspace to develop a Koopman operator~\cite{brunton2016koopman}. In their work, they showed that the state matrix found by approximating the Koopman operater for specific systems can be used to generate control laws with linear quadratic regression (LQR). 
Further more, Proctor et al.~\cite{proctor2014dynamic} expanded on their previous work on DMD to create a DMD method with control to extract low-order control models from higher-dimensional systems building the base DMD algorithm. This DMD variant with control, called DMDc, is demonstrated to show positive results in the analysis of infections disease data. 
Proctor et al.~\cite{proctor2016generalizing} introduced a method of finding the Koopman operator of a system that takes into account system's inputs and control (KIC) based on prior work on DMDc.
Work has also been done to improve the accuracy of the Koopman operater by 'lifting' the states of the system to a set of observables such that a solution can be found where the data required from DMD is limited, extending the DMD algorithm (EDMD) \cite{Williams_2015}.
Korda et al.~\cite{Korda_2018} presented the use of lifting nonlinear dynamics on augmented states with EDMD and explored finding control laws for systems using model predictive control (MPC).
There have also been work done to optimize the methods in which the lifting of dynamics is done through the optimization of the dictionary of functions used to lift the state variables through EDMD with dictionary learning~\cite{Li_2017}. 

Applications of the Koopman operator have been explored in the field of robotics, where the operator is used to develop closed-loop controllers for pendulum systems~\cite{Abraham_2017}. 
Kaiser et al.~\cite{kaiser2020datadriven} explores the application of the Koopman operator theory in generating energy-based control using a DMD variant with control, extended dynamic mode decomposition with control, and compared it against another method of finding the Koopman operator called the Koopman reduced order nonlinear identification and control (KRONIC).
Conversly, the Koopman operator is not a panacea to solving dynamical systems; a recent work by Gonzalez et al.~\cite{antikoopmanism} aimed at demonstrating some of the weakness of the Koopman operator and the DMD algorithms and provides alternate approaches to reconstruction.

This paper will illustrate the mathematics behind the Koopman operator theory and show how augmenting the states of a dynamical system can create accurate linear models through this process. 
It then shows how the DMD algorithm can be used to find a linear approximation of dynamical systems before exploring variants of DMD to create more accurate approximations of the Koopman operator when representing system dynamics and other methods of using DMD to generate control of a system represented via the Koopman operator.  
Numerical simulations of two classical dynamical systems, the inverted pendulum and cart-pole systems, are then studied using the DMD algorithms. 
% The remainder of this paper provides analysis and discussion on the system reconstruction using the variants of the DMD algorithm to approximate the Koopman operator. 

\iffalse
\begin{figure}[t]
    \centering
    \includegraphics[width=0.95\linewidth]{figs/koopmanvis.pdf}
    \caption{The Koopman operator 'lifts' the system's dynamics from the state space to a higher dimensional augmented state where the observed dynamics are both linear and near infinite. }
    \label{fig:kop_ill}
\end{figure}

 \fi
\section{Preliminaries}

%\subsection{Continuous and Discrete Systems}

This work pertains to the implementation of Koopman operator theory on dynamical systems to propagate non-linear dynamical equations with a data-driven approximation method. Conventionally, a non-linear dynamical system consists of a set of states and a function or rule that governs how the states propagate either forward in time or with respect to each other~\cite{DataDriven,DMDbook}. This can be described in a continuous and a discrete method. 

For a generic continuous system
\begin{align}
    \frac{d}{dt}\mathbf{x}(t) &= F(\mathbf{x}(t),t;\mu),\label{eq:dynamics}
\end{align}
$\mathbf{x}(t) \in \mathbb{R}^n$ is the vector holding the states of the dynamical system at time $t$, $n$ is the number of states that define the system, $\mu$ is a set of parameters for the system dynamics, and $F(\cdot)$ is the rule describing the evolution of the state in a continuous sense. 
These continuous-time dynamics can also be modeled in a respective discrete-time representation where the system can be evaluated at every finite time interval, $\Delta t$, which can otherwise be viewed as $\mathbf{x}_k = x(k\Delta t)$ with the subscript $k$. 
The evolution of a dynamic system in a discrete-time flow map can be formally portrayed as
\begin{align}
    \mathbf{x}_{k+1} &= f(\mathbf{x}_k),    \label{eq:2}
\end{align}
by collecting the states at time $t_k$, $k = 1,2,...,m$, for $m$ time steps, where $x_k$ is a $n$-dimensional column vector of system states and $x_{k+1}$ is the states of the system in the time step following $x_k$~\cite{proctor2016generalizing,proctor2014dynamic}.

 These rules governing how the states of the system advance tend to be nonlinear equations to best simulate most of the systems in practice. An unfortunate quality of nonlinear systems is that their dynamics are difficult to solve analytically. As a result, modern control practices tend to turn to approximations in order to produce a high-fidelity controller for the system.
 However, if a dynamical system can be expressed with by a linear rule, we can then achieve more accurate predictions of how a system advances in time. We intend to show this using the Koopman operator theory.

\subsection{Koopman Operator}
The Koopman operator is a linear but infinite dimensional operator that can be defined for an autonomous, discrete time, dynamical system. Unlike the function $f$ shown earlier in \eqref{eq:2}, the Koopman operator advances the measurements of the changes in dynamics over time\cite{DataDriven,proctor2016generalizing,Korda_2018}. Consider the evolution of a dynamical system $\mathbf{x}_{k+1} = f(\mathbf{x}_k)$ where the rule $f$ maps the state space onto itself, i.e., $f: \mathbb{R}^n \rightarrow \mathbb{R}^n$; using Koopman operator theory we define a different rule, $g: \mathbb{R}^n \rightarrow \mathbb{M}^{n_y}$, where $n_y$ is the dimension of a nearly infinite column vector defining the dimension of the observable of $\mathbf{x}$ at a given time step. This implies that the Koopman operator is defined for all observables, meaning that the Koopman operator $\mathcal{K}$ is also infinite-dimensional. We will discuss how we can make this a more reasonable attribute later on in this discussion. Here, $g$ is a real-valued, scalar, measurement function which is an element of an infinite-dimensional Hilbert space called an observable. The Koopman operator acts on this observable such that
\begin{align}
\mathcal{K}_t g =& g \circ F(\mathbf{x}(t)),   \label{eq:3}\\
\mathcal{K}_{\Delta t}g(\mathbf{x}_k) =& g(f(\mathbf{x}_k)).        \label{eq:4}
\end{align}
In \eqref{eq:3}, the Koopman operator evolves the observable with respect to time in the continuous fashion while in \eqref{eq:4} the Koopman operator evolves the discrete-time dynamics with respect to $\Delta t$. More details on the connections between these two representations can be found in \cite{DataDriven,proctor2016generalizing,Korda_2018,kaiser2020datadriven,kooopman_sys_ctrl}. Equations \eqref{eq:3} and \eqref{eq:4} allow us to define an analogue for continuous and discrete-time dynamical systems respectively as
\begin{align}
    \frac{d}{dt}g =& \mathcal{K}g \label{eq:5},\\
    g(\mathbf{x}_{k+1})=&\mathcal{K}_{\Delta t}g(\mathbf{x}_k) \label{eq:6}.
\end{align}

This work explores the application of the Koopman operator on discret-time dynamical systems. Due to the linear nature of the Koopman operator, we can perform an eigen decomposition of $\mathcal{K}$ such that
\begin{align}
    \mathcal{K}\varphi_j(\mathbf{x}_k) = \lambda_j\varphi_j(\mathbf{x}_k), \label{eq:7}
\end{align}
where $\lambda$ and $\varphi$ are the Koopman eigenvalue and eigenvector describing the evolution of the Koopman operator. Considering this, the observable, $g(x)$, can be expanded as
\begin{align}
    g(\mathbf{x}) &=
    \begin{bmatrix}
    g_1(\mathbf{x})\\g_2(\mathbf{x})\\\vdots\\g_i(\mathbf{x})
    \end{bmatrix}
    &= \sum_{j=1}^{\infty}\varphi_j(\mathbf{x})\bf{v}_j , \label{eq:8}
\end{align}
where $\bf{v}$ is a coefficient called the Koopman mode associated with its corresponding Koopman eigenvector. This allows us to consider Koopman modes as a projection of the observable:
\begin{equation}
    \bf{v}_j =
    \begin{bmatrix}
    \langle \varphi_j,g_1\rangle\\\langle \varphi_j,g_2\rangle\\\vdots\\\langle \varphi_j,g_i\rangle
    \end{bmatrix}\label{eq:9}.
\end{equation}
By combining \eqref{eq:7} and \eqref{eq:8}, one can show the relationship between the observable, $g$, and the Koopman operator, $\mathcal{K}$. However, for all practical proposes, using an infinite-dimensional vector and operator is not always feasible so an approximation of the Koopman operator, $\mathcal{K}$, will be typically used and found from collected data on the system. This approximation allows us to combine the terms in  \eqref{eq:7} and \eqref{eq:8} to show the relationship between the Koopman modes in $g$ and the Koopman eigenvalues in $\mathcal{K}$
\begin{equation}
    \mathcal{K}g(\mathbf{x}_k) \approx \mathcal{K}g(\mathbf{x}_k) = g(f(\mathbf{x}_k)) = \sum_{j=1}^{\infty}\lambda_j\varphi_j(\mathbf{x})\bf{v}_j.
    \label{eqn:koopman_operator_show}
\end{equation}
This shows that the Koopman operator is an iterative set of triples, $\lambda_j, \varphi_j$ and $\bf{v}_j$, all of which make up the Koopman mode decomposition\cite{DataDriven,proctor2016generalizing,Korda_2018,kaiser2020datadriven,Mezic_2005,Mezic2013}. 

\begin{comment}
I want to use this section someewhere but I am not sure how to best put it in because I think it is a good high-level summary of what the koopman operator does before I introduce the mathematics

This changes how the original dynamical system represented by (\mathbf{x},k,$f$) is transformed into a new dynamical system defined by ($F$,k,$K$). In a graphical sense instead of mapping a system's state from $\mathbf{x}$ to $f(\mathbf{x})$ in the same plane as if this were a simple linear transformation, the Koopman operator updates the observables of the original dynamical system which would take the plane that the original state resides, move it to the new observable plane which transforms the the observable to a new value and then brings that value to a corresponding location on the original state plane. %(An image would be good here like in the Rowley paper)
. The evolution between observables is the Koopman operator. \\
To perform the Koopman analysis of a given system, the Koopman eigenvalues($\mu_k$) and eigenfunctions ($\psi_k$) which capture the long term dynamics of the observables and the Koopman modes ($\nu_k$) which are vectors that allows us to reconstruct the system's state as a linear combination of the Koopman eigenfunctions.
\end{comment}

\subsection{Dynamic Mode Decomposition (DMD)}
DMD is a data-driven method of exploring the behavior of complex systems. Using measurements from numerical simulations or in experiments, we can attempt to find the most dominant dynamic characteristics of the system. DMD does this by finding the dynamic modes or eigen modes and eigenvalues of the proposed system. \\
DMD acts on the assumption that the state of a system is connected to the next by
\begin{equation} 
    \mathbf{x}_{k+1} = \mathbf{A}\mathbf{x}_k, \label{eq:linprop}
\end{equation}
where $\mathbf{x}(t) \in \mathbb{R}^n$ and $A \in \mathbb{R}^{n\times n}$ and is the matrix describing the evolution of the state in a continuous-time manner. Simulated or experimental measurements for $\mathbf{x}_k$ are then collected at regular time intervals of $\Delta t$ to become snapshots to be used in a discrete time system. These snapshots are collected and stored in sequence like the following:
\begin{align}
\bf{\mathbf{X}} =& \begin{bmatrix}
\vline && \vline && && \vline\\
\mathbf{x}_1 && \mathbf{x}_2 && \ldots && \mathbf{x}_{m-1}\\
\vline && \vline && && \vline
\end{bmatrix},\label{eq:data1}\\
\bf{\mathbf{X}'} =& \begin{bmatrix}
\vline && \vline && && \vline\\
\mathbf{x}_2 && \mathbf{x}_3 && \ldots && \mathbf{x}_{m   }\\
\vline && \vline && && \vline
\end{bmatrix}\label{eq:data2},
\end{align}
where $\bf{\mathbf{X}'}$ is the time-shifted snapshot of matrix $\bf{\mathbf{X}}$ such that
\begin{align}
    \bf{\mathbf{X}'} \approx A\bf{\mathbf{X}} \label{eq:13}.
\end{align}
There is no set number of snapshots required for DMD but a sufficiently large number is needed for the application and is closely related to the approximation of the Koopman operator. For DMD we are attempting to find the $\mathbf{A}$ in \eqref{eq:13} using the snapshots of the system's states. $\mathbf{A}$ can be approximated by
\begin{align}
    \mathbf{A} =& \mathbf{X}'\mathbf{X}^+\label{eq:XXpmt},
\end{align}
where $^+$ is the Moore-Penrose pseudoinverse. Considering that $A$ has a large $n$ that normal calculation would be a large computational load we can perform Singular Value Decomposition (SVD) on the snapshots to find the dominant characteristics of the pseudoinverse of $\bf{\mathbf{X}}$
\begin{align}
    \bf{\mathbf{X}} \approx U \Sigma V^*\label{eq:svd},
\end{align}
where $\mathbf{U} \in \mathbb{R}^{n\times r}$, $\mathbf{\Sigma} \in \mathbb{R}^{r\times r}$ and $\mathbf{V} \in \mathbb{R}^{n\times r}$ where $*$ denotes the conjugate transpose. $r$ is the reduced rank of the SVD approximation of $\bf{\mathbf{X}}$. From SVD we can rearrange the singular values and the eigenvectors of $\bf{\mathbf{X}}$ with the forward snapshot to  $A$ that satisfies \eqref{eq:13}:
\begin{align}
    % A\approx \Bar{A} =& \bf{\mathbf{X}}'\tilde{V}\tilde{\Sigma}^{-1}\tilde{U}^*\label{eq:17}'
     \mathbf{A}\approx \Bar{\mathbf{A}} =& \bf{\mathbf{X}}'\tilde{V}\tilde{\Sigma}^{-1}\tilde{U}^*\label{eq:17}.
\end{align}
%where $\Bar{A}$ is an approximation of $A$ seen in \eqref{eq:svd}.
However in practice since the dimensionality of $A$ is so large we have to approximate it as
\begin{align}
    \tilde{\mathbf{A}} = \mathbf{U}^*\mathbf{X}'V\mathbf{\Sigma}^{-1}\label{eq:appxA},
\end{align}
where $\tilde{A}$ is a $r$ rank linear model of the dynamical system such that
\begin{align}
    \tilde{\mathbf{x}}_{k+1} =& \tilde{\mathbf{A}}\tilde{\mathbf{x}}_k.
\end{align}
Here $\tilde{\mathbf{x}} : \in\mathbb{R}^{r\times 1}$ is also a reduced order of the state $\mathbf{x}$ which can be reconstructed by
\begin{align}
    \tilde{\mathbf{x}}:\mathbf{x}=\tilde{\mathbf{U}}\tilde{\mathbf{x}}.
\end{align}
The next step is to find the spectral decomposition of $\tilde{\mathbf{A}}$:
\begin{align}
    \tilde{\mathbf{A}}\mathbf{W}=\mathbf{W}\mathbf{\Lambda}\label{eq:eigenequal},
\end{align}
where $\mathbf{\Lambda}$ and $\mathbf{W}$ are the corresponding eigenvalues (DMD modes) and eigenvectors of the full rank system dynamics matrix allowing us to recover the full state system dynamics in a computationally efficient manner since $\tilde{\mathbf{A}} \in\mathbb{R}^{r\times r}$ where $r<<n$.
We can now reconstruct the eigendecomposition in \eqref{eq:svd} for $\mathbf{A}$ from $\mathbf{W}$ and $\mathbf{\Lambda}$ with the corresponding eigen vectors given by the columns of $\mathbf{\Phi}$ such that
\begin{align}
    \mathbf{\Phi} =& \mathbf{X}'\mathbf{V\Sigma}^{-1}\mathbf{W}.
\end{align}
From $\mathbf{\Phi}$  we can no reconstruct our approximation of the time dynamics of $\mathbf{x}(t)$ by projecting our approximations into a future solution
\begin{align}
    \mathbf{x}(t) \approx \sum_{k=1}^{r}\boldsymbol{\phi}_k\exp{(\omega_kt)}b_k = \mathbf{\Phi}\exp{(\mathbf{\Omega} t)}b, \label{eq:eigen_reconstruction}
\end{align}
where $b_k$ is the initial amplitude of each DMD mode, $\boldsymbol{\phi}$ is the columns that make up $\mathbf{\Phi}$ and $\mathbf{\Omega}$ is the diagonal matrix of the eigenvalues $\omega$ in $\omega_k = \ln{(\lambda_k)}/\Delta t$\cite{DMDbook,DataDriven,kooopman_sys_ctrl,proctor2014dynamic,H_Tu_2014}.

\subsection {Extended Dynamic Mode Decomposition}

Extended DMD is almost the same algorithm as the standard DMD one, however the method of which we deploy EDMD is by using the observables of the system to create a dictionary to pass through the normal DMD algorithm.\cite{Williams_2015,DataDriven,DMDbook,williams2015kernelbased}
By performing regression on this new augmented vector containing linear and non-linear measurements we can make a different approximation of the original system dynamics.
This augmented vector is constructed in the following manner
\begin{equation}
    \mathbf{y} = \mathbf{\Theta}^T(\mathbf{x}) = \begin{bmatrix}
\boldsymbol{\theta}_1(\mathbf{x})\\ 
\boldsymbol{\theta}_2(\mathbf{x})\\ 
\vdots \\
\boldsymbol{\theta}_p(\mathbf{x})
\end{bmatrix},\label{eq:augref}
\end{equation}
where $p$ is the rank of the augmented state such that $p>>n$. Here $\mathbf{\Theta}$ is the collection of measurements of the system possibly containing the original state of the system, $\mathbf{x}$, as well as nonlinear measurements. Once $\mathbf{y}$ is found, two data matrices are created in the same manner seen above in the DMD algorithm \eqref{eq:data1} and \eqref{eq:data2}. From here a best-fit linear operator $\mathbf{A}_Y$ is found that maps \eqref{eq:data1} onto \eqref{eq:data2} 
\begin{equation}
    \mathbf{A}_Y =\operatorname*{argmin}_{\mathbf{A}_Y} \left \|{\mathbf{Y}'-\mathbf{A}_Y \mathbf{Y}}  \right \|=\bold{Y'Y^+}.
\end{equation}
This regression can than be written in terms of the original data matrices $\mathbf{\Theta}^T(\mathbf{x})$ as
\begin{equation}
    \mathbf{A}_Y =\operatorname*{argmin}_{\mathbf{A}_Y} \left \|\mathbf{\Theta}^T(\mathbf{x}')-\mathbf{A}_Y \mathbf{\Theta}^T(\mathbf{x})  \right \|=\mathbf{\Theta}^T(\mathbf{x}')(\mathbf{\Theta}^T(\mathbf{x}))^+\label{eq:redaugdata}.
\end{equation}
$\mathbf{A}_Y$ is then the basis upon which we can derive the Koopman operator. However, $\mathbf{\Theta}$ may not necessarily span the same subspace as the Koopman operator and may consist of different eigenvalues and eigenvectors of the Koopman operator when is why verification and re-validation techniques need to be used to show that the EDMD model is properly fit to the actual system.

\subsection{Approximating Koopman Eigenfunctions from Data}
Beginning with the observable matrices in \eqref{eq:augref}, the Koopman eigenfunction can be approximated as
\begin{equation}
    \varphi (\mathbf{x}) \approx \sum_{k=1}^p \boldsymbol{\theta}_k(\mathbf{x})\boldsymbol{\xi}_k = \mathbf{\Theta}(\mathbf{x})\boldsymbol{\xi}\label{eq:modexpan},
\end{equation}
where by using \eqref{eq:augref}, we can find
\begin{equation}
    \begin{bmatrix}
\lambda\varphi(\mathbf{x}_1)\\ 
\lambda\varphi(\mathbf{x}_2)\\ 
\vdots \\
\lambda\varphi(\mathbf{x}_m)\\ 

\end{bmatrix} = \begin{bmatrix}
\varphi(\mathbf{x}_2)\\ 
\varphi(\mathbf{x}_3)\\ 
\vdots\\
\varphi(\mathbf{x}_{m+1}) 

\end{bmatrix}\label{eq:eigprop}.
\end{equation}
We start seeing a connection between the Koopman operator and the process we followed to approximate a linear evolution seen in DMD. We can expand \eqref{eq:modexpan} as
\begin{align}
    [\lambda\mathbf{\Theta}(\mathbf{X})-\mathbf{\Theta}(\mathbf{X}')]\boldsymbol{\xi} = 0\label{eq:boop}.
\end{align}
Afterwards, we reduce \eqref{eq:boop} using a best least-squares fit to get
\begin{align}
    \lambda\boldsymbol{\xi} =\mathbf{\Theta}(\mathbf{X})^+\mathbf{\Theta}(\mathbf{X}')\boldsymbol{\xi}\label{eq:lsfaugstate}.
\end{align}
Now we can compare \eqref{eq:redaugdata} and \eqref{eq:lsfaugstate}; \eqref{eq:redaugdata} is the transpose of the latter so that the left eigenvectors are now the right eigenvectors. Comparing to \eqref{eq:XXpmt} we can use the eigenvectors $\boldsymbol{\xi}$ of $\mathbf{\Theta}^+\mathbf{\Theta}'$ to find the coefficients of the eigenfunction $\varphi(\mathbf{x})$ that represents the basis of $\mathbf{\Theta}(\mathbf{x})$. We now can confirm that the predicted eigenfunctions behave linearly by comparing them to the predicted dynamics in \eqref{eq:eigprop}. We do this to determine if the regression done in \eqref{eq:lsfaugstate} yields proper eigenvalues and eigenvectors that span the Koopman invariant subspace for the system\cite{Williams_2015}.

\subsection{Methodology}
The dynamical system model and the methods we used to discover the Koopman operator are presented in the section. As discussed in the previous section, the method that will be used to determine a finite-dimensional approximation of the Koopman operator will be through DMD and EDMD. This approach has been used in ~\cite{Brunton_2016,proctor2016generalizing,kaiser2020datadriven}. As stated in the previous section, DMD computes the transformation matrix $\mathbf{A}$ that shows the evolution of two measurement snapshot pairs by finding the corresponding eigenvectors and eigenvalues that satisfy
\begin{align}
    \mathbf{A}\mathbf{v}_j =& \lambda_j\mathbf{v}_j.
\end{align}
Now assuming that $\mathbf{A}$ contains a complete set of eigenvectors, each measurement can be found by the eigenvectors of $\mathbf{A}$ like so:
\begin{align}
    g(x_k) = \sum_{j=1}^{n}c_{jk}\mathbf{v}_j.
    \label{eqn:findmode}
\end{align}

In the case that these snapshot or augmented pairs in the case of EDMD and DMD, the DMD modes and eigenvalues of \eqref{eqn:findmode} correspond to the Koopman modes of \eqref{eqn:koopman_operator_show}. Since the principle concept of the Koopman operator, in the realm of data-driven identification, is to advance the states of the system by advancing the dynamics linearly.

\section{Results and Analysis}

In this work we will numerically evaluate controlled and uncontrolled systems to collect sufficiently large sets of data, convert it into snapshot pairs to run through the DMD and EDMD algorithms and use the decomposition to reconstruct linear models of the data to compare against the simulated systems.

\subsection{Sample Systems}
The classical systems that were used in this work to approximate the Koopman operator are the inverted pendulum and cart pole system~\cite{DataDriven}. These systems are used due to their periodic and transient behaviors and dynamics are well known. 

\subsubsection{Inverted Pendulum}
\begin{figure}
    \centering
    \includegraphics[trim={150 180 150 120}, clip, width=0.7\linewidth]{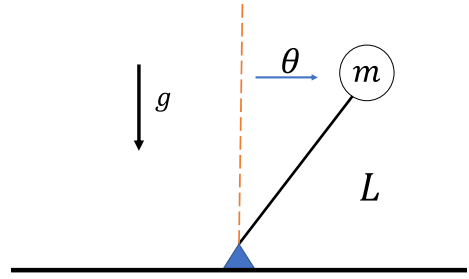}
    \caption{ Classical inverted pendulum system with mass $m$ = 1 kg and arm length of $L$ = 2 m. In this model the mass is free to swing about the pivot 360 degrees and will be treated as an undamped system. }
    \label{fig:pend}
\end{figure}

In the inverted pendulum system seen in figure \ref{fig:pend} we will start with an initial state of $[\pi/4,0]$ and in the controlled examples we attempted to right the pendulum mass to a terminal state of $[0,0]$. In this simulation we only controlled the $\dot{\theta}$ term with a performance index of $[0,10]$ and a control cost of 1 for our LQR gain. 

\subsubsection{Cart-Pole}
\begin{figure}
    \centering
    \includegraphics[trim={200, 180, 100, 100}, clip, width=0.7\linewidth]{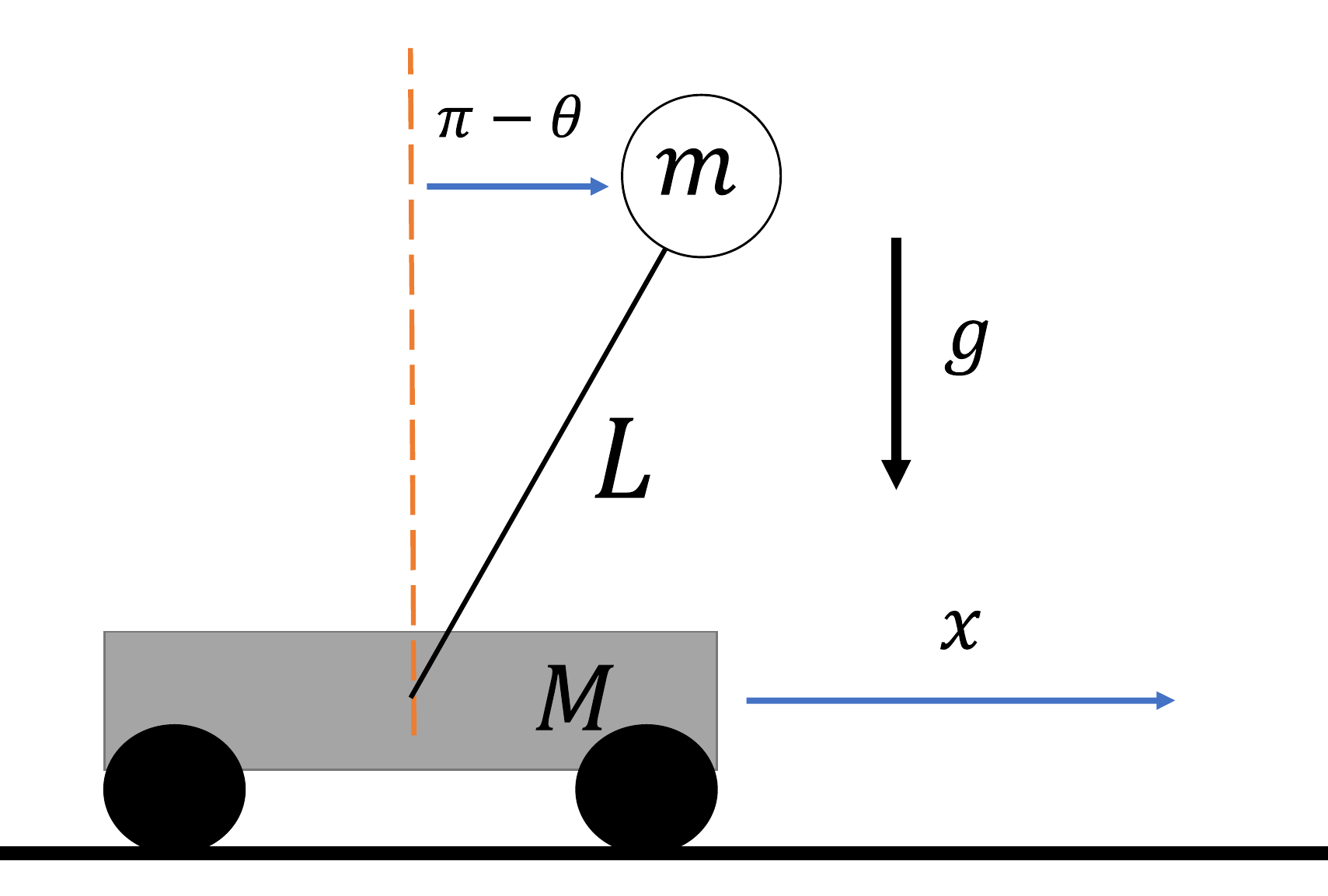}
    \caption{ Classical cart pole system with pendulum mass $m$ = 1kg, cart mass $M$ = 5 kg, and arm length $L$ = 2 m. Like the pendulum system in Fig.~\ref{fig:pend}, the end mass is free to swing about the pivot 360 degrees and both the cart and arm pivot are friction-less and undamped. }
    \label{fig:cartpole}
\end{figure}
 In the cart-pole system seen in Fig.~\ref{fig:cartpole} we started with an initial state of $[-1,0,\pi,0]$ and attempted to right the cart-pole system after a brief translation to $[1,0,\pi,0]$. In this simulation we created a control law using LQR with a performance index of $[5, 10, 0, 0]$ and a control cost of 1. 

\subsection{Implementation}
\begin{figure}
    \centering
    \includegraphics[trim={0, 40, 0, 50}, clip, width=0.95\linewidth]{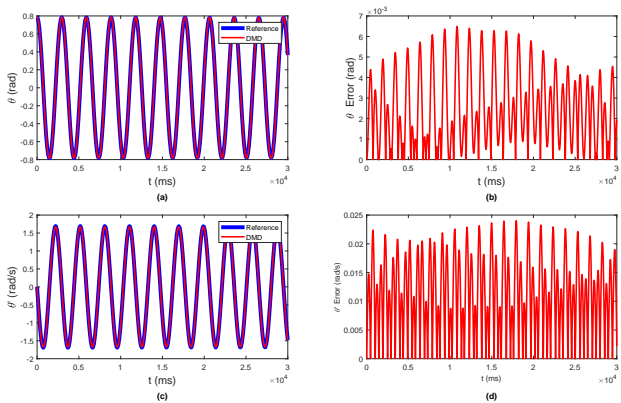}
    \caption{DMD performed on a pendulum system with uncontrolled data.}
    \label{fig:DMD_States_no_control}
\end{figure}
\begin{figure}
    \centering
    \includegraphics[trim={0, 40, 0, 50}, clip, width=0.95\linewidth]{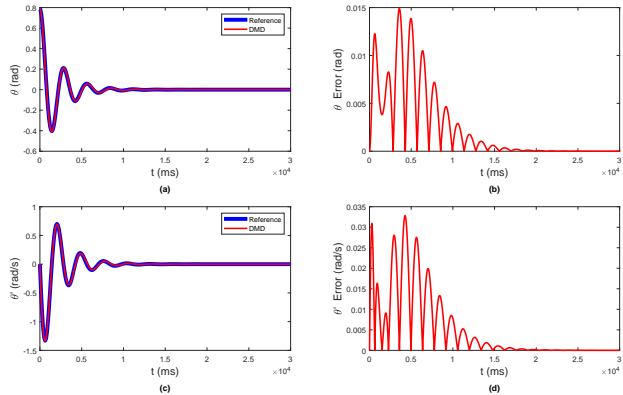}
    \caption{DMD performed on a pendulum system with with LQR controlled data.}
    \label{fig:DMD_States_w_control}
\end{figure}
We can see in Figs.~\ref{fig:DMD_States_no_control} and \ref{fig:DMD_States_w_control} that the reconstruction of the original state space succeeded with minimal errors. The periodic nature of the pendulum system can still be seen in the error analysis and it recedes in time. This is most likely due to the immediate transient behavior of the reference data, though the over all behavior appears to be accurately represented. We can see consistent low error in the reconstructed controlled and uncontrolled pendulum data. 
\begin{figure}
    \centering
    \includegraphics[trim={0, 40, 0, 50}, clip, width=0.95\linewidth]{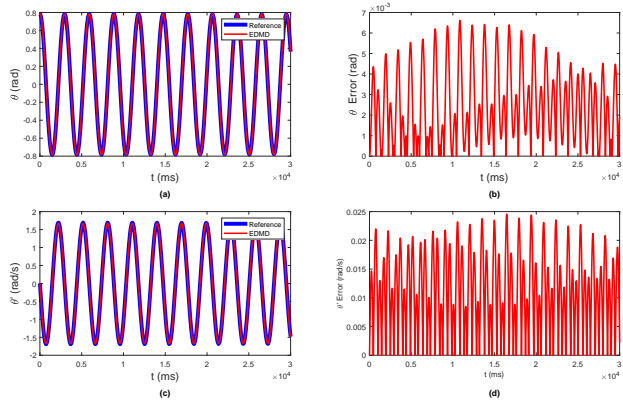}
    \caption{EDMD performed on a pendulum system with uncontrolled data lifted by a second-order polynomial basis.}
    \label{fig:EDMD_States_no_control}
\end{figure}
\begin{figure}
    \centering
    \includegraphics[trim={0, 40, 0, 50}, clip, width=0.95\linewidth]{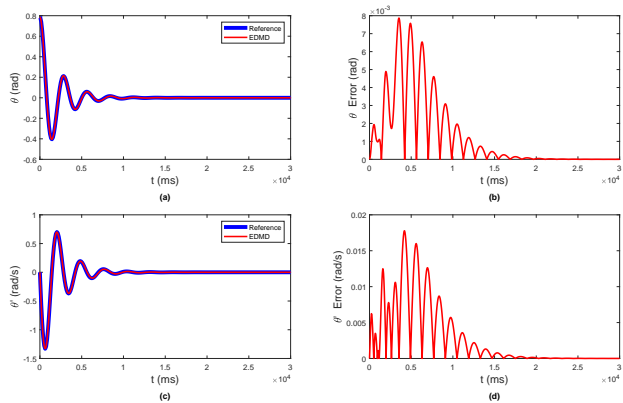}
    \caption{EDMD performed on a pendulum system with with LQR controlled data lifted by a second-order polynomial basis.}
    \label{fig:EDMD_States_w_control}
\end{figure}
In Figs.~\ref{fig:EDMD_States_no_control} and \ref{fig:EDMD_States_w_control} there is the repeated behavior of the DMD reconstruction seen in Figs.~\ref{fig:DMD_States_no_control} and \ref{fig:DMD_States_w_control}. In both sets of the EDMD reconstruction we can see similar or reduced error in reconstruction. 

\begin{figure}
    \centering
    \includegraphics[width=0.95\linewidth]{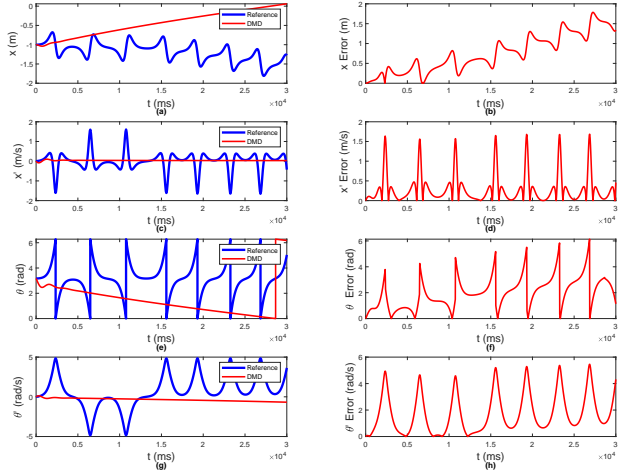}
    \caption{DMD performed on a cart-pole system with uncontrolled data.}
    \label{fig:cart_pole_DMD_States_no_control}
\end{figure}
\begin{figure}
    \centering
    \includegraphics[width=0.95\linewidth]{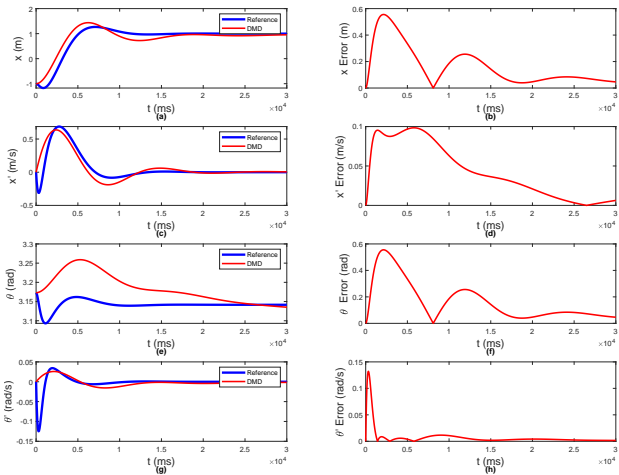}
    \caption{DMD performed on a cart-pole system with with LQR controlled data.}
    \label{fig:cart_pole_DMD_States_w_control}
\end{figure}
In Figs.~\ref{fig:cart_pole_DMD_States_no_control} and \ref{fig:cart_pole_DMD_States_w_control} we can see significant losses in reconstructing the behavior of the uncontrolled cart-pole system. Figure ~\ref{fig:cart_pole_DMD_States_w_control} shows that some of the characteristics of the reference data have been captured in the reconstruction as the reference data in this process is smoother and has less periodic and transient characteristics. 

\subsection{Discussions}

\begin{figure}
    \centering
    \includegraphics[trim={0, 40, 0, 50}, clip, width=0.95\linewidth]{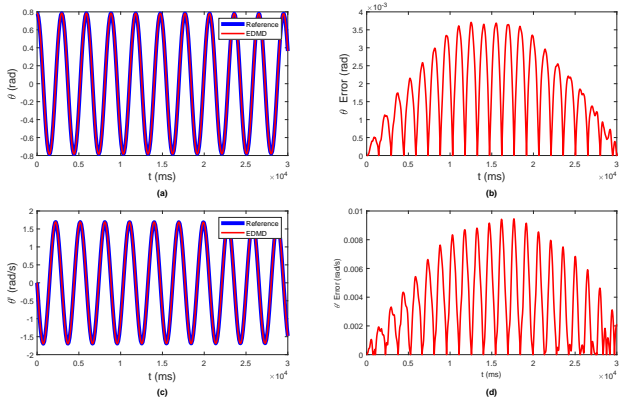}
    \caption{EDMD performed on a pendulum system with uncontrolled data using a first-order Fourier basis to lift the state data. }
    \label{fig:EDMD_fourier_States_no_control}
\end{figure}

\begin{figure}
    \centering
    \includegraphics[trim={0, 40, 0, 50}, clip, width=0.95\linewidth]{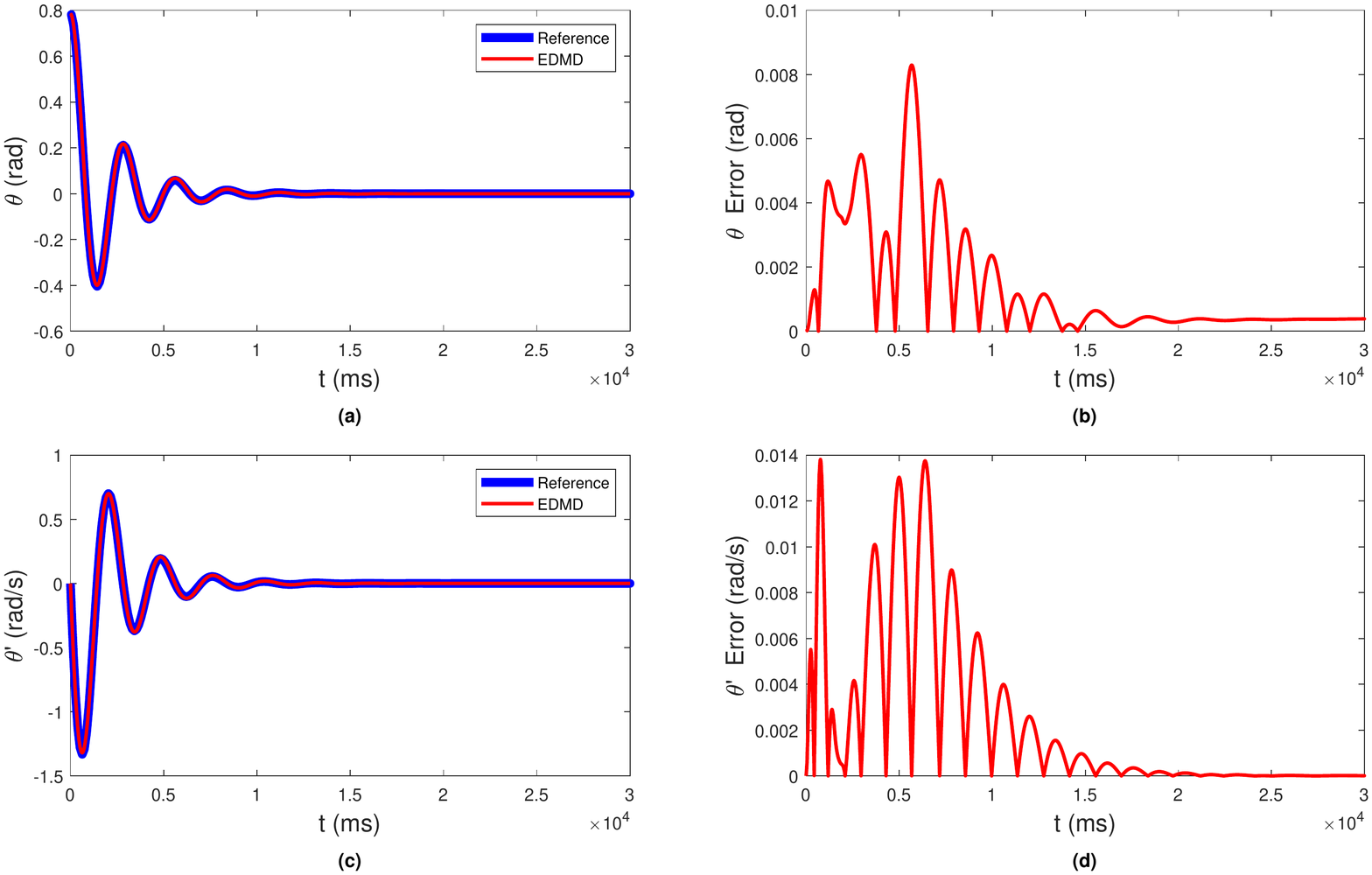}
    \caption{EDMD performed on a pendulum system with LQR controlled data using a first-order Fourier basis to lift the state data. }
    \label{fig:EDMD_fourier_States_w_control}
\end{figure}

\begin{figure}
    \centering
    \includegraphics[trim={0, 40, 0, 50}, clip, width=0.95\linewidth]{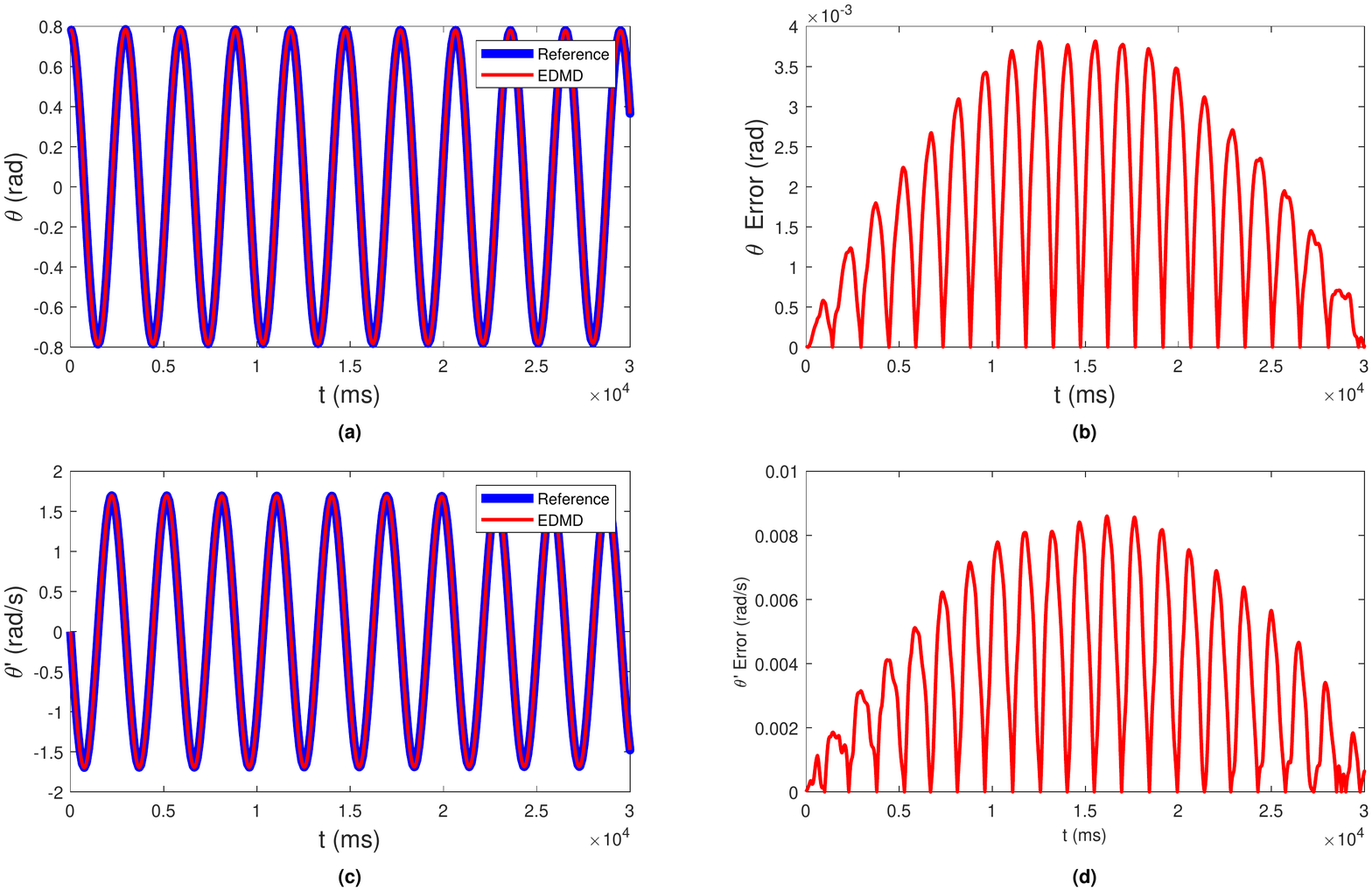}
    \caption{EDMD performed on a pendulum system with uncontrolled data using a second-order Fourier basis to lift the state data. }
    \label{fig:EDMD_fourier2_States_no_control}
\end{figure}

We can see in the reconstruction results that the inverted pendulum system is a good candidate for the system reconstruction using the approximated Koopman operator. We can see a marginal improvement in the reconstruction of the state data of the controlled data when we improve the DMD algorithm to EDMD thought the improvement yet no change in the uncontrolled data. This shows that the periodic behavior is not improved by a polynomial lifting function. Additionally when moving to a higher order polynomial (or Fourier; see Figs.~\ref{fig:EDMD_fourier_States_no_control} and \ref{fig:EDMD_fourier_States_w_control}) we can additional improvements in the reconstruction of the original data when the full lifted states are used. Oddly enough, as seen in Fig.~\ref{fig:EDMD_fourier2_States_no_control}, we see no further improvement in the overall fidelity of the reconstruction of the state space using a second-order Fourier fitting.

\begin{figure}
    \centering
    \includegraphics[trim={0, 40, 0, 50}, clip, width=0.95\linewidth]{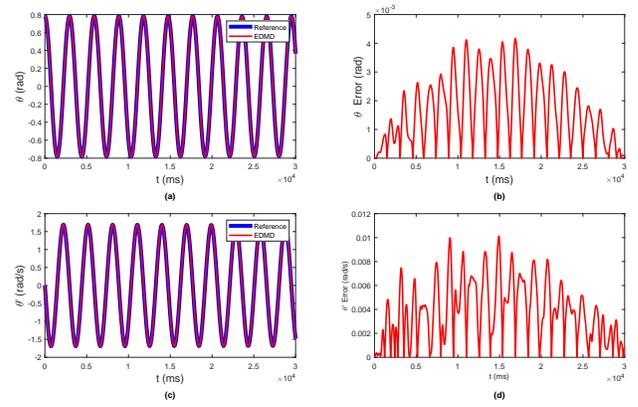}
    \caption{EDMD performed on a pendulum system with uncontrolled data lifted by a third-order polynomial basis with no truncation.}
    \label{fig:EDMD2_States_no_control}
\end{figure}
\begin{figure}
    \centering
    \includegraphics[trim={0, 40, 0, 50}, clip, width=0.95\linewidth]{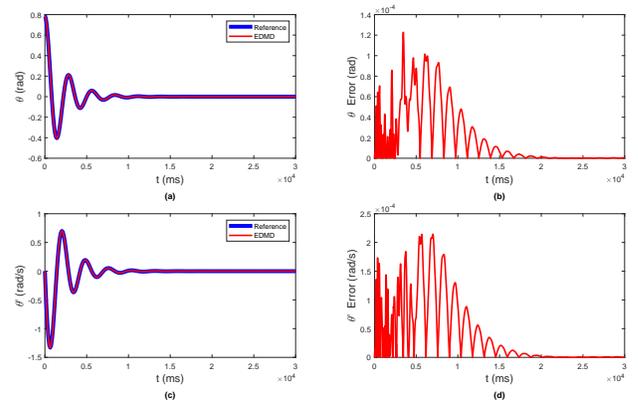}
    \caption{EDMD performed on a pendulum system with with LQR controlled data lifted by a third-order polynomial basis with no truncation.}
    \label{fig:EDMD2_States_w_control}
\end{figure}

\begin{figure}
    \centering
    \includegraphics[width=0.95\linewidth]{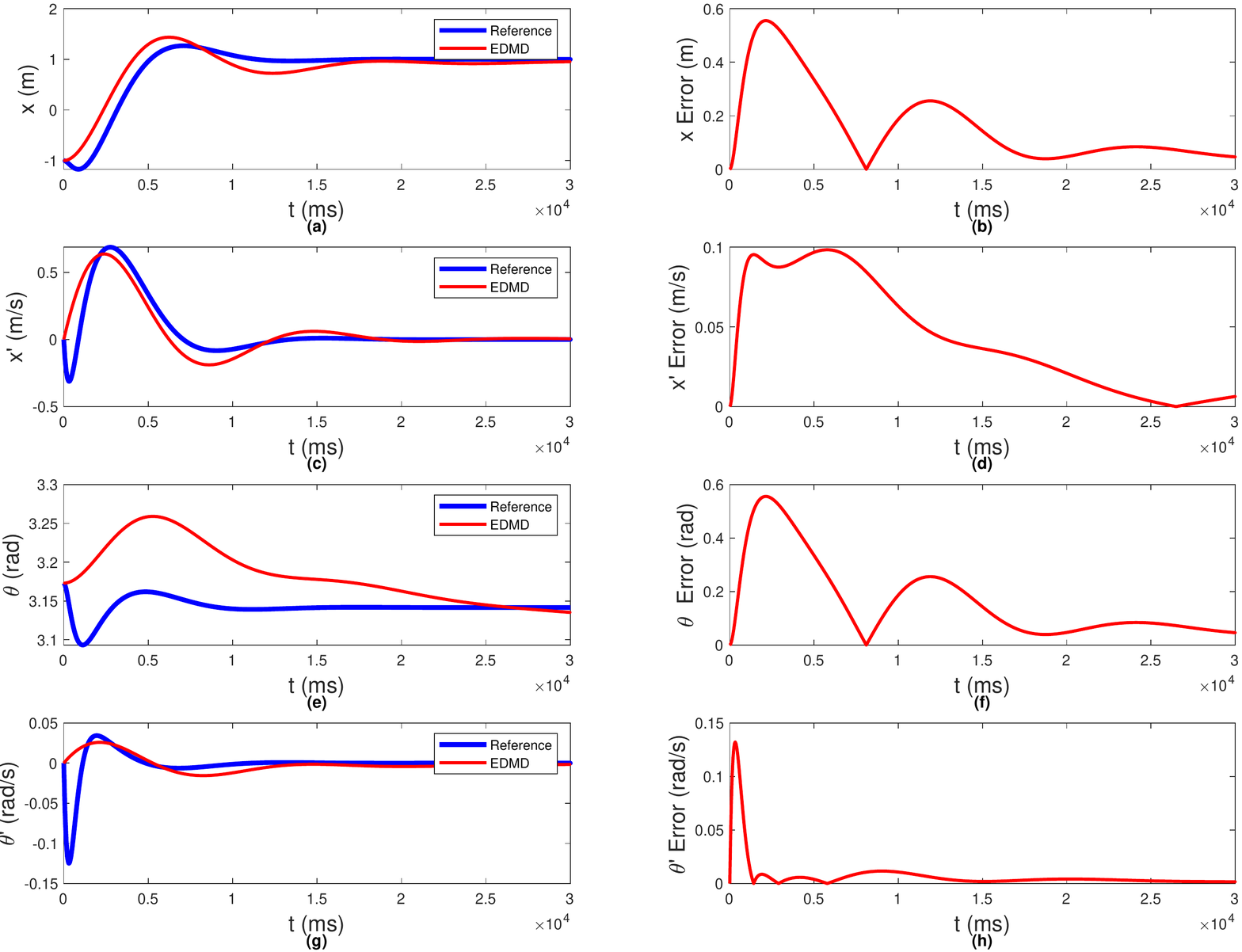}
    \caption{EDMD performed on a cart-pole system with LQR controlled data using a second-order polynomial basis to lift the state data. }
    \label{fig:cart_pole_EDMD_States_w_control}
\end{figure}

\begin{figure}
    \centering
    \includegraphics[width=0.95\linewidth]{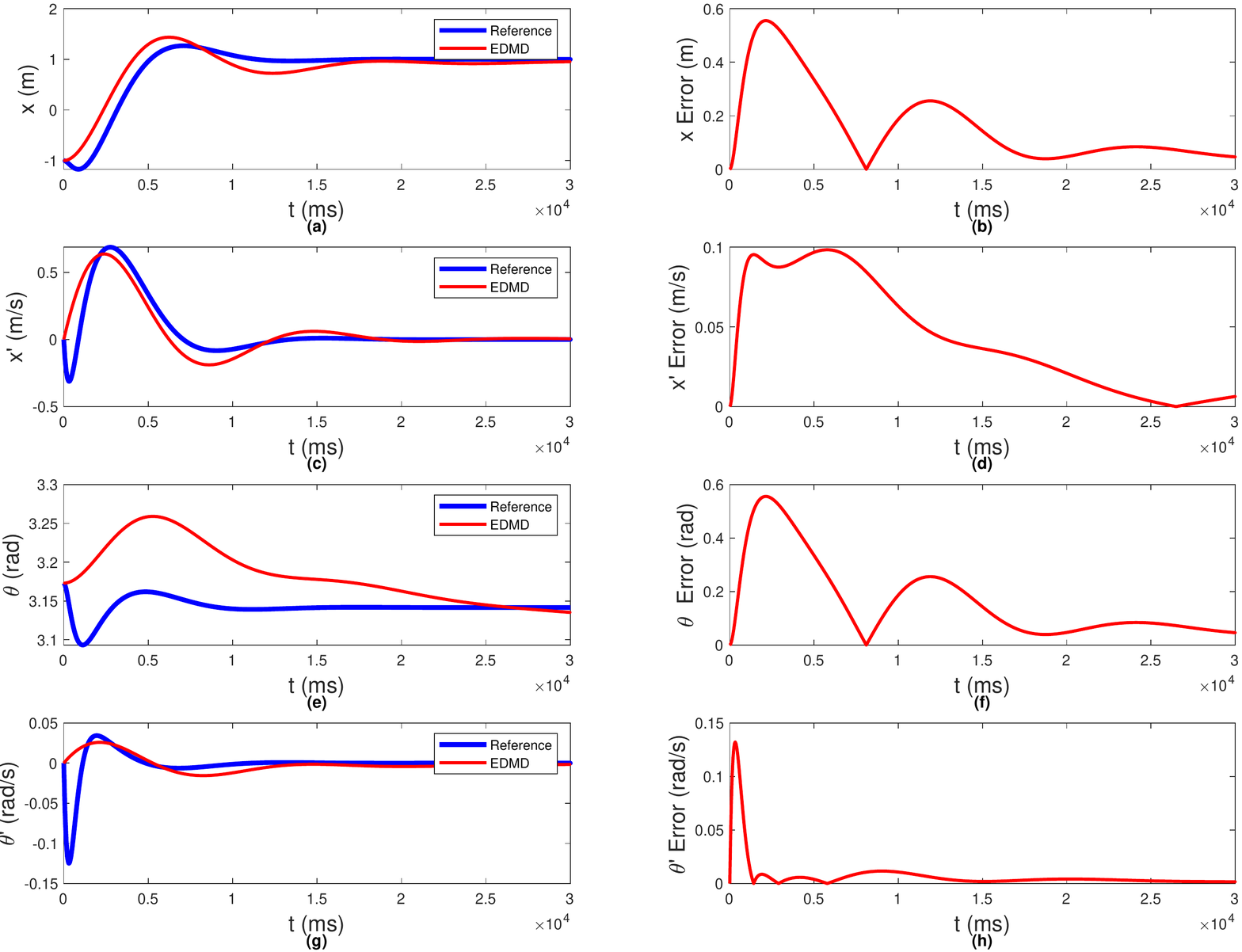}
    \caption{EDMD performed on a cartpole system with LQR controlled data using a first-order Fourier basis to lift the state data. }
    \label{fig:cart_pole_EDMD_fourier_States_w_control}
\end{figure}

Moving on to the second system in this work, we can see significant errors in the reconstruction process as we are adding more transient components into the system dynamics with the extra degrees of freedom of the cart-pole system. As shown in Fig.~\ref{fig:cart_pole_DMD_States_no_control}, DMD (and by extension, EDMD) did not capture the uncontrolled behavior of the system in question however when the strong periodic behavior is removed in the controlled data set we find an improved reconstruction of the system however weaker that that shown in the inverted pendulum system. Additionally we note that the initial response generated by LQR appears to be missed by the reconstruction possibly showing that more gradual control may yield better results in the reconstruction. Additionally we can see no further improvement on the reconstruction of the system when we attempted to use the EDMD algorithm in Figs.~\ref{fig:cart_pole_EDMD_States_w_control} and \ref{fig:cart_pole_EDMD_fourier_States_w_control} using either a polynomial or Fourier basis to lift the original state space. 
%~\ref{fig:cart_pole_EDMD_States_no_control},

%~\ref{fig:cart_pole_EDMD_fourier_States_no_control},

\begin{figure}
    \centering
    \includegraphics[trim={0, 40, 0, 50}, clip, width=0.95\linewidth]{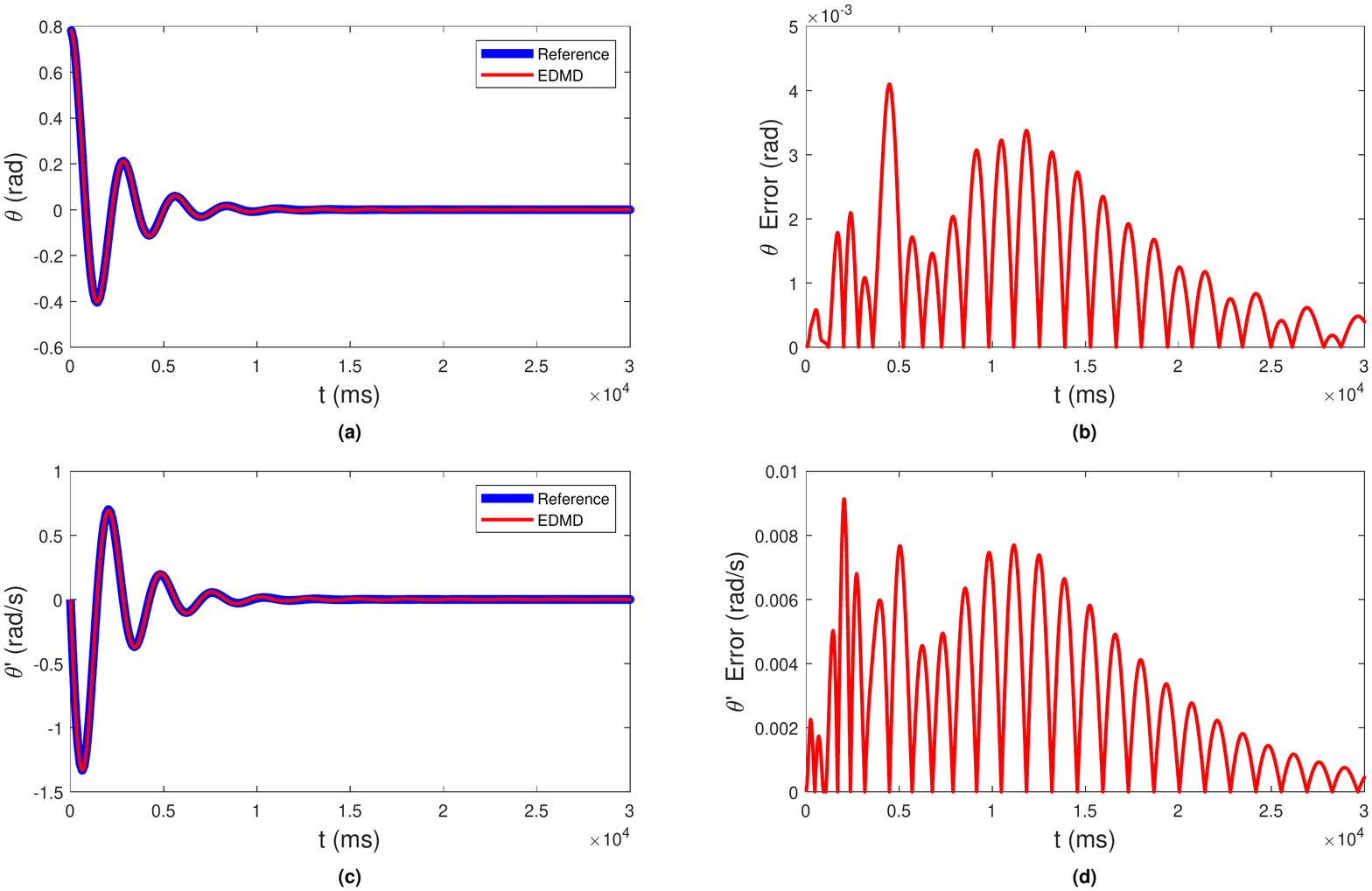}
    \caption{EDMD performed on a pendulum system with LQR controlled data using a second-order Fourier basis to lift the state data. }
    \label{fig:EDMD_fourier2_States_w_control}
\end{figure}


\section{Conclusions}

This work presented an introduction on the usage of the Koopman operator theory approach to model simple, nonlinear dynamical systems with a linear operator.
In this process, we simulated the states of classical non-linear systems and feed that data as a set of snapshots through the DMD algorithm, and its variants, to reconstruct the state data using the Koopman operator. 
Through the usage of these algorithms on an inverted pendulum and cart-pole system we can approximate the Koopman operator to recreate high fidelity linear approximations of these systems and provide increased accuracy by augmenting the states of the original data.
These methods showed weakness in reconstructing systems that act with coupled states that behave in different matters as seen in the cart-pole system. 
In future work, we intend to explore alternate lifting functions to represent the state spaces of our systems and explore alternative methods of finding the Koopman operator of systems with multiple coupled states to improve performance to develop linear control laws.

\iffalse
\section*{Acknowledgment}
The authors thank the Office of Naval Research and the National Science Foundation for funding this project.
\fi

%
 \bibliographystyle{IEEEtran}
\bibliography{ranbib/Gregory}
\iffalse
\appendices
\section*{Appendix} \label{Appendix}

\begin{figure}
    \centering
    \includegraphics[width=0.95\linewidth]{koopman-vtol/figs/DMDc_States_no_control.pdf}
    \caption{DMDc performed on a pendulum system with uncontrolled data.}
    \label{fig:DMDc_States_no_control}
\end{figure}
\begin{figure}
    \centering
    \includegraphics[width=0.95\linewidth]{koopman-vtol/figs/DMDc_States_w_control.pdf}
    \caption{DMDc performed on a pendulum system with with LQR controlled data.}
    \label{fig:DMDc_States_w_control}
\end{figure}

\begin{figure}
    \centering
    \includegraphics[width=0.95\linewidth]{koopman-vtol/figs/cart_pole_EDMD_States_no_control.pdf}
    \caption{EDMD performed on a cart-pole system with uncontrolled data using a second-order polynomial basis to lift the state data. }
    \label{fig:cart_pole_EDMD_States_no_control}
\end{figure}

\begin{figure}
    \centering
    \includegraphics[width=0.95\linewidth]{koopman-vtol/figs/cart_pole_EDMD_fourier_States_no_control.pdf}
    \caption{EDMD performed on a cart-pole system with uncontrolled data using a first-order Fourier basis to lift the state data. }
    \label{fig:cart_pole_EDMD_fourier_States_no_control}
\end{figure}
\fi


%\ifCLASSOPTIONcaptionsoff
%  \newpage
%\fi

\end{document}